\documentclass[12pt]{article}
\usepackage{latexsym}
\begin{document}
\baselineskip 5mm
\title{Ricci Curvature, Diameter and Fundamental Groups}
\author{Wen-Haw Chen$^*$ and Jyh-Yang Wu\thanks{%
Partially supported by Taiwan NSC grants} }
\date{ }
\maketitle
\footnotetext{E-Mail: whchen@mail.thu.edu.tw;
                      jywu@math.ccu.edu.tw}

\begin{abstract}
    In this note we discuss the fundamental groups and diameters of positively Ricci curved $n$-manifolds. We use a method combining the results about equivarient Hausdorff convergence developed by Fukaya and Yamaguchi with the Ricci version of splitting theorem by Cheeger and Colding to give new information on the topology of compact manifolds with positive Ricci curvature. Moreover, we also obtain a weak Margulis's lemma for manifolds under a lower Ricci curvature bound. 
\end{abstract}
   \section{Introduction}    
     This paper concerns about the obstruction problems for compact manifolds with positive Ricci curvature. In dimension 2, this problem is easy to understand since only the projective plan $\mbox{\boldmath$RP$}^2$ and the 2-sphere admit metrics with positive curvature. In dimension 3, Hamilton showed in {\bf [H]} that a compact $3$-manifold with positive Ricci curvature also admits a metric with a constant sectional curvature of $+1$, and then is covered by the $3$-sphere. In general, a classical result of Myers shows that the fundamental group of a compact positively Ricci curved manifold must be finite. Moreover, since $SU(n)$ has positive Ricci curvature, any finite group will occur as the fundamental group of some manifold with positive Ricci curvature.

     The problem 5 listed in Lecture Series 4 in [{\bf [P2]}, p.105] states: {\em Consider a compact positively Ricci curved manifold, what can be said about the fundamental group depending only on the dimension $n$, except for the fact that it is finite?} In this paper, we give for this problem a partial answer, which was conjectured by the second-named author in {\bf [W]}. On the other hand, we show that the diameter of the universal Riemannian covering space of a compact positively Ricci curved manifold $M$ can not be too larger than that of $M$. The second-named author also conjectured this property for positively curved manifold.\\
\\
{\bf Theorem A.} {\em Given $n\geq 2$, there exist constants $p_n$ and $C_n$ depending only on $n$ such that if a                                                compact Riemannian $n$-manifold $M^n$ has the Ricci curvature $Ric_{M^n}>0$,  then}\\
                 (a) {\em the first betti number $b_1(M^n,\mbox{\boldmath$Z$}_p)$ with p-cyclic group
                     coefficient  $\mbox{\boldmath$Z$}_p$ satisfies $b_1(M^n,\mbox{\boldmath$Z$}_p)\leq n-1$ for all prime                                           $p\geq p_n$}, and\\
                 (b) the ratio of diameters satisfies
\[\frac{diam(\tilde{M}^n)}{diam(M^n)}\;\;<\;\; C_n,\] {\em where $\tilde{M}^n$ is the  universal covering of $M^n$.} \\
\\
{\bf Remark 1.1.} In view of the flat $n$-torus $T^{n}$, one has $b_1(T^{n},\mbox{\boldmath$Z$}_p)=n$ for all prime $p$ and the canonical Euclidean $n$-space $\mbox{\boldmath$R$}^n$ is its universal covering space. Hence our assumption for curvature is optimal. In fact, Fukaya and Yamaguchi give in [Corollary 0.9 in {\bf [FY1]}] that if a compact Riemannian $n$-manifold $M$ with sectional curvature $K_{M}$ and diameter $diam(M)$ satiefies $K_{M}diam(M)^2>-\epsilon_{n}$ for some constant $\epsilon_{n}$ depending only on $n$, then $b_1(M^n,\mbox{\boldmath$Z$}_p)\leq n$ for all $p\geq p(n)$ and the maximal case $b_1(M^n,\mbox{\boldmath$Z$}_p)=n$ occurs only when $M^{n}$ is diffeomorphic to a torus. Also they obtained in [Corollary 0.11 in {\bf [FY1]}] that $diam(\tilde{M})/diam(M)$ is uniformly bounded by a constant depending only on $n$ provided the fundamental group $\pi_1(M)$ is additionally finite. Theorem A extends their results to manifolds with positive Ricci curvature. Note that if $n\leq 3$ and $Ric_M>0$, then it is covered by spheres as discussed above. Hence our Theorem A holds for these manifolds.\\
\\
{\bf Remark 1.2.} Here is an application of Theorem A. It is well-known that every finite group $G$ can be the fundamental group of a compact $4$-manifold. If we take $G=S_m$ to be the permutation group of $m$ elements and consider the $4$-manifolds $M_m^4$ with fundamental group $S_m$, then $M_m^4$ admits no metric with positive Ricci curvature for large $m$.\\

 From the argument in the proof of Theorem A, we have the following weak Margulis's lemma under a lower Ricci curvature bound. Recall that the {\em length of polycyclicity} of a solvable group $G$ is the smallest integer $m$ for which $G$ admits a  filtration
 \[\{e\}=G_m\subset G_{m-1}\subset\ldots\subset G_1\subset G_0=G\]
 such that each $G_i/G_{i-1}$ is cyclic.\\
\\
{\bf Theorem B (A weak Margulis's Lemma).} {\em There exists a positive number $\delta_n$ depending only on $n$ and satisfying the following: Let $(M^n,p)$ be a complete pointed Riemannian $n$-manifold with $Ric_{M^n}\ge -(n-1)$. Then there exists a point $p'\in B_p(1/2)$ such that the image of the inclusion homomorphism
\[\Gamma'\;=\;Im[\pi_1(B_{p'}(\delta_n))\,\to\,\pi_1(B_p(1))]\]
admits a subgroup $\Lambda'\subset\Gamma'$ with}\\
(1) {\em $\;\;[\Gamma':\Lambda']<w_n$, where $w_n$ depends only on $n$};\\
(2) {\em $\;\;\Lambda'$ is solvable with length of polycyclicity $\leq n$}.

{\em In particular, if a complete Riemannian $n$-manifold $M$ has $Ric_Mdiam(M)^2> -(n-1)\delta_n$, then $\pi_1(M)$ is almost solvable. That is, $\pi_1(M)$ contains a solvable subgroup of finite index.}\\
\\ 
{\bf Remark 1.4.} Gromov conjectured in {\bf [G]} that there is a positive number $\epsilon_n$ depending only on $n$ such that if a compact Riemannian $n$-manifold with almost nonnegative Ricci curvature $Ric_Mdiam(M)^2>-\epsilon_n$, then $\pi_1(M)$ is {\em almost nilpotent}. Fukaya and Yamaguchi showed in {\bf [FY1]} that Gromov's conjecture is true under the condition $K_Mdiam(M)^2>-\epsilon_n$. By taking a solvable subgroup in place of a nilpotent subgroup, they gave a generalized Margulis's lemma in [Theorem A2.1 in {\bf [FY1]}]. In {\bf [FY1]}, it was suggested that in order to extend the Margulis's lemma to Ricci case, it "only" need to establish splitting theorem and volume convergence theorem under almost nonnegative Ricci curvature bound. However, this is not enough to extend Fukaya and Yamaguchi's result to Ricci case by using the original arguments since their argument depends mainly on the existence of a fibration with the property of {\em almost Riemannian submersion}, and it is not true in general that one can construct such a fibration under a lower Ricci curvature bound. 
Therefore, the techniques of Fukaya and Yamaguchi cannot carry over directly to the manifolds with a lower Ricci 
curvature bound. 

In our approach, we need to work harder on the induction steps and prove a Technical lemma 3.1 in section 3, which is weaker 
than the original one in [Theorem 7.1 in {\bf [FY1]}]. Indeed, we are not able to obtain the Margulis' lemma under a lower Ricci 
curvature bound. We can only obtain a weaker version of the Margulis' lemma for only " one point " in the Riemannian manifold 
under consideration. This is sufficient for us to prove the solvability theorem for almost nonnegatively Ricci curved manifolds. 
Though we are still unable to obtain the nilpotency result, Theorem B confirms, in some sense, the almost solvability version of 
Gromov's conjecture.\\

The remainder of this paper is divided into four sections: In section 2, we mention the main tools including the theory of pointed equivarient convergence and the splitting theorems for our proof of Theorem A and Theorem B. Especially, we establish Corollary 2.6 by combining the above two tools. In section 3, we prove a Technical lemma, which extends the solvability theorem in [Theorem 7.1 in {\bf [FY1]}]. In section 4, we give a proof of Theorem B by using the Technical lemma. Then Theorem A can be proved in section 5 by Theorem B. \\
   \section{Equivariant Pointed Hausdorff Convergence and the Splitting Theorem}     
We first recall the notion about the  equivariant pointed Hausdorff convergence in [3]. Let ${\cal M}$ be the set of all isometry classes of pointed metric spaces $(X,p)$ such that, for each $D>0$, the ball $B_p(D)$ around $p$ with radius $D$ is relatively compact and such that $X$ is a length space. Denote ${\cal M}_{eq}$ the set of triples $(X, \Gamma, p)$, where $(X,p)\in {\cal M}$ and $\Gamma$ is a closed subgroup of isometries of $X$. Put 
$\Gamma(D)\;=\;\{\gamma\in\Gamma\mid d(\gamma p,p)<D\}$.\\
\\
{\bf Definition 2.1.} Let $(X,\Gamma,p)$, $(Y,G,q)\in {\cal M}_{eq}$. An $\varepsilon$-{\it equivariant pointed Hausdorff approximation} is a triple $(f,\phi,\psi)$ of maps $f:B_p(1/\varepsilon)\to Y$, $\phi:\Gamma(1/\varepsilon)\to G(1/\varepsilon)$ and $\psi:G(1/\varepsilon)\to\Gamma(1/\varepsilon)$ such that\\
(2.1.1) $\;\;\;f(p)=q$;\\
(2.1.2) $\;\;\;$the $\varepsilon$-neighborhood of $f(B_p(1/\varepsilon))$ contains $B_q(1/\varepsilon)$;\\
(2.1.3) $\;\;\;$if $x\,\;y\in B_p(1/\varepsilon)$, then $\mid d(f(x),f(y))-d(x,y)\mid <\varepsilon$;\\
(2.1.4) $\;\;\;$if $\gamma\in\Gamma(1/\varepsilon)$, $x\in B_p(1/\varepsilon)$, $\gamma x\in B_p(1/\varepsilon)$, then \[d(f(\gamma x),\phi(\gamma)(f(x))<\varepsilon;\]
(2.1.5) $\;\;\;$if $\lambda\in G(1/\varepsilon)$, $x\in B_p(1/\varepsilon)$, $\psi(\mu)(x)\in B_p(1/\varepsilon)$, then \[d(f(\psi(\mu)(x)),\mu f(x))<\varepsilon.\]
Hereafter the notion $\lim_{i\to\infty}(X_i,G_i,x_i)\;=\;(Y,G,y)$ means
\[\lim_{i\to\infty}d_{eH}((X_i,G_i,x_i),(Y,G,y))\;=\;0\]
, where $d_{eH}$ denotes the {\em the equivariant pointed Hausdorff distance}. For the sake of brief, we also denote $d_{eH}$ to be $d_H$ in the remainder of this paper.\\

The following theorem comes from [Proposition 3.6 in {\bf [FY1]}].\\
\\
{\bf Theorem 2.2.} {\em Let $(X_i,\Gamma_i,p_i)\in {\cal M}_{eq}$, $(Y,q)\in {\cal M}$. Suppose that 
$\lim_{i\to\infty}(X_i,p_i)\;=\;(Y,q).$
Then $G$ and a subsequence $k_i$ can be found such that $(Y,G,q)\in{\cal M}_{eq}$ and
$\lim_{i\to\infty}(X_{k_i},\Gamma_{k_i},p_{k_i})\;=\;(Y,G,q).$}\\

The following theorem is shown in [Theorem 4.2 in {\bf [FY2]}] by Fukaya and Yamaguchi. Its proof can be found in [Appendix A.1 in {\bf [FY1]}].\\
\\
{\bf Theorem 2.3.} {\em Let $(X_i,\Gamma_i,p_i)$, $(Y,G,q)\in {\cal M}_{eq}$ be such that 
$\lim_{i\to\infty}(X_i,\Gamma_i,p_i)\;=\;(Y,G,q),$
and $G'$ be a normal subgroup of $G$. Assume that}\\
(2.3.1) {\em $\;\;\;G/G'$ is discrete.}\\
(2.3.2) {\em $\;\;\;Y/G$ is compact.}\\
(2.3.3) {\em $\;\;\;\Gamma_i$ is discrete and free and $X_i$ is simply connected.}\\
(2.3.4) {\em $\;\;\;G'$ is generated by $G'(R_0)$ for some $R_0>0$}.\\
{\em Then there exists a sequence of normal subgroups $\Gamma_i$ of $\Gamma$ such that} \\
(2.3.5) {\em $\;\;\;\lim_{i\to\infty}(X_i,\Gamma_i',p_i)\;=\;(Y,G',q).$}\\
(2.3.6) {\em $\;\;\;\Gamma_i/\Gamma_i'$ is isometric to $G/G'$ for sufficiently large $i$.}\\
(2.3.7) {\em $\;\;\;G/G'$ is finitely presented.}\\
(2.3.8) {\em $\;\;\;\Gamma_i'$ is generated by $\Gamma_i'(R_0+\varepsilon_i)$ for some $\varepsilon_i\to 0$}\\

The next result due to Cheeger and Colding shows in [Theorem 6.64 in {\bf [CC1]}] that the limit space of a sequence of complete pointed-Riemannian $n$-manifolds with almost nonnegative Ricci curvature will split provided that it contains a line. \\
\\
{\bf Theorem 2.4.} {\em Let $(M_i^n,p_i)$ be a sequence of complete pointed-Riemannian $n$-manifolds. Denote $B_{p_i}(R_i)$ be the open $R_i$-ball in $M_i^n$ around $p_i$ and $R_i\to\infty$ as $i\to\infty$. Let $(X,p_{\infty})\in {\cal M}$ with $\lim_{i\to\infty}(B_{p_i}(R_i),p_i)\;=\;(X,p_{\infty})$. Suppose $Ric_{B_{p_i}(R_i)}\ge -\epsilon_i^2$, where $\epsilon_i\to 0$ as $i\to\infty$, and $X$ contains a line. Then $X$ splits, isometrically, $X\,=\,\mbox{\boldmath$R$}\times X'.$}\\

Combine Theorem 2.2 with Theorem 2.4, one can extend straightforward, as in [Corollary 5.3  and Theorem 5.4 in {\bf [FY1]}], to Corollary 2.5 and Corollary 2.6 respectively. Corollary 2.6 is especially important for the proof of our main result.\\
\\
{\bf Corollary 2.5.} {\em Let $(M_i,p_i)$, $B_{p_i}(R_i)$ and $(X,p_{\infty})$ be as in Theorem 2.4. Suppose $G_i$ is a closed subgroup of $Isom(B_{p_i}(R_i))$ such that $diam(B_{p_i}(R_i)/G_i)\leq D$ for some constant $D$. Then there exists a subsequence $k_i$ that
  \[\lim_{i\to\infty}(B_{p_{k_i}}(R_{k_i}),G_{k_i},p_{k_i})\,=\,(\mbox{\boldmath$R$}^{\ell}\times Y, G, p_{\infty})\]                
, where $Y$ is a compact metric space and $G$ is a closed subgroup of $Isom(\mbox{\boldmath$R$}^{\ell}\times Y)$ with $\ell\leq n$.}\\
\\
{\bf Corollary 2.6.} {\em Let $(M_i,p_i)$ be a sequence of complete pointed-Riemannian $n$-manifolds with $Ric_{M_i}\ge -(n-1)$. Suppose  $\lim_{i\to\infty}(M_i,p_i)\;=\;(X,p_{\infty})$, where $(X,p_{\infty})\in {\cal M}$. Then for every $x\in X$ there exists sequences $y_i\in X$, $q_i\in M_i$ and $r_i\to\infty$ as $i\to\infty$ such that}\\
   (2.6.1) {\em $\;\;\;y_i\to x$, $q_i\to x$ as $i\to\infty$},\\
   (2.6.2) {\em $\;\;\;\lim_{i\to\infty}((X,r_id_X),y_i)=(\mbox{\boldmath$R$}^k,can,0)=\lim_{i\to\infty}((M_i,r_ig_i),q_i)$},\\
   (2.6.3) {\em $\;\;\;k\leq n$,\\
    where $d_X$ and $g_i$ are the original metric of $X$ and $M_i$ respectively.}\\
\\
{\bf Remark 2.7} As in the proof of [Theorem 5.4 in {\bf [FY1]}], Corollary 2.6 can be proved by blowing up the metrics at most finite times and using Theorem 2.4. Note that for a given convergent sequence $\delta_i\to 0$ one can always find a sequence $r_i\to\infty$ in Corollary 2.6 such that $r_i\delta_i\to 0$.\\

Let $Y$ be a compact metric space and $G$ a closed subgroup of $Isom(\mbox{\boldmath$R$}^{\ell}\times Y)$. Since $G$ preserves the splitting $\mbox{\boldmath$R$}^{\ell}\times Y$, the projection $\phi:G\to Isom(\mbox{\boldmath$R$}^{\ell})$ is well defined. The following theorem was shown in [Lemma 6.1 in {\bf [FY1]}].\\
\\
{\bf Theorem 2.8.} {\em For each $\varepsilon>0$ there exists a normal subgroup $G_{\varepsilon}$ of $G$ such that}\\
(2.8.1) {\em $\;\;\;G/G_{\varepsilon}$ is discrete;}\\
(2.8.2) {\em $\;\;\;$there exists an exact sequence
        $1\to G_{\varepsilon}\to G\to\Lambda\to 1,$\\
        where $\Lambda$ contains a finite-index free abelian subgroup of rank not greater than $\dim(\mbox{\boldmath$R$}^{\ell}/\phi(G))$;}\\
(2.8.3) {\em $\;\;\;$for every $g\in G_{\varepsilon}$ and every $x\in\mbox{\boldmath$R$}^{\ell}\times Y$ there exists $g_1,\ldots,g_s\in G_{\varepsilon}$ satisfying}\\
$\;\;\;(i)$ {\em $g=g_s\ldots g_1$,}\\
$\;\;\;(ii)$ {\em $d(g_ig_{i-1}\ldots g_1(x),g_{i-1}\ldots g_1(x))<\varepsilon$ for all $1\leq i\leq s$.}\\

The group $G_{\varepsilon}$ was constructed in {\bf [FY1]} as follows: Let $K=Ker(\phi)$, which acts on $Y$. Set 
$\hat{K}_{\varepsilon}\;=\;\{g\in K\mid d(g(x),x)<\varepsilon\,\forall\,x\in Y\}.$
Let $K_{\varepsilon}$ be the group generated by $\hat{K}_{\varepsilon}$. Since $K_{\varepsilon}$ is normal in $G$, the natural projection $\pi:G\to G/K_{\varepsilon}$ is defined. Define 
$G_{\varepsilon}\;\,=\;\,\pi^{-1}((G/K_{\varepsilon})_0)$
, where $(G/K_{\varepsilon})_0$ denotes the identity component of $G/K_{\varepsilon}$.\\
\\ 
{\bf Remark 2.9.} If the limit space $Y$ is an Alexandrov space, then Fukaya and Yamaguchi showed in {\bf [FY2]} that $Isom(Y)$ is in fact a Lie group. Thus $G$ is a Lie group and one can take $G_0$ as $G_{\varepsilon}$ for every $\varepsilon$. On the other hand, Cheeger and Colding announced a result that if $Ric_{M_i^n}\ge -(n-1)$ and $Vol(B_{p_i}(1))\ge v>0$ for all $i$ and all $p_i\in M_i$, then the isometry group of the limit space is a Lie group. Note that the construction of $G_{\epsilon}$ in Theorem 2.8 is independent of curvature and volume.\\
   \section{A Technical Lemma}                           

  The following Technical lemma plays a very important role in our approach to prove Theorem A and Theorem B, and it has its own interest for investigating  manifolds with lower Ricci curvature bounds. It can be viewed as a {\em weak} Ricci version of [Theorem 7.1 in {\bf [FY1]}].\\
\\ 
{\bf Technical Lemma 3.1.} {\em For given positive integers $n$ and $k$, $n\geq k$ and a positive number $\mu_0$, there exists positive numbers $\epsilon=\epsilon_{n,k}(\mu_0)$, $w=w_{n,k}$ and a function $\tau(\epsilon)=\tau_{n,k,\mu_0}(\epsilon)$ with $\lim_{\epsilon\to 0}\tau(\epsilon)=0$ such that if $(M^n,p)$ and $(N^k,q)$ are pointed-Riemannian manifolds of dimension $n$ and $k$ respectively such that}\\
(3.1.1) {\em $\;\;\;Ric_M\geq -(n-1)$, $Ric_N\geq -(n-1)$ and $inj(N)>\mu_0>0$,}\\
(3.1.2) {\em $\;\;\;d_{GH}((M,p),(N,q))<\epsilon$,\\
 where $d_{GH}$ denotes the Gromov-Hausdorff distance,
   then there exists a map $f:M\to N$ with $f(p)=q$ satisfying the following:}\\
(3.1.3) {\em $f$ is a continuous $\tau(\epsilon)$-Hausdorff approximation such that $f_*:\pi_1(M,p)\to\;\pi_1(N,q)$ is          surjective;}\\
(3.1.4) {\em Let $V=B_q(\frac{\mu_0}{2})$ be the ball around $q$ with radius $\frac{\mu_0}{2}$. Set                                             $U=f^{-1}(V)$. Then there is a normal subgroup $H$ of the fundamental group $\Gamma=\pi_1(U)$ of $U$                                              such that}\\ 
 (i) {\em $\;\;\;H$ is a solvable subgroup of $\Gamma$ with length of polycyclicity $\leq n-k$,}\\
  (ii) {\em $\;\;\;[\Gamma : H]\leq w_{n,k}$.}\\
\\
{\bf Remark 3.2.} Recently, Sormani and Wei considered in {\bf [SW]} the group $\bar{\pi}_1(Y)$ of deck transforms of the                   universal cover $Y$ of the Gromov-Hausdorff limit of compact manifolds $\{M_i^n\}$ with $Ric_{M_i^n}\ge                                                      (n-1)H$ and $diam(M_i^n)\leq D$ for some $H\in \mbox{\boldmath$R$}$ and $D>0$. They showed that for $n_0$                   sufficient large depending on $Y$, there is a surjective homeomorphism                                                     $\Phi_i:\pi_1(M_i)\to\bar{\pi}_1(Y)$, for $i\ge n_0$. In Technical Lemma 3.1, we consider Gromov-Hausdorff convergence to LGC($\rho$)-space to obtain a similar result.

We divided the proof of Technical lemma 3.1 into the following two parts. The first part shows the existence of the map $f$ satisfying (3.1.3).\\
\\
{\bf Proof of (3.1.3).} The construction of such a Hausdorff approximation $f$ depends heavily on the assumption that $inj_{q}(N)\ge\mu_0$. We say a function   $\rho:[0,r)\to [0,\infty)$ to be a contractibility function provided:(i) $\rho(0)=0$, (ii) $\rho(\epsilon)\ge\epsilon$,   (iii) $\rho(\epsilon)\to 0$ as $\epsilon\to 0$, (iv) $\rho$ is non-decreasing. Then a metric space $X$ is said to be an {\em LGC$(\rho)$-space with a contractibility function $\rho$} if for every $\epsilon\in [0,r]$ and $x\in X$ the ball $B_x(\epsilon)$ is contractible inside $B_x(\rho(\epsilon))$. Since $inj_q(N)>\mu_0>0$, $(N,q)$ is an LGC($\rho$)-space with contractibility function $\rho(s)=s$ defined in $[0,\mu_0/2]$.

 Choose $\epsilon$ be such that $8(n+3)^{2}\epsilon<\mu_0$. Since $d_{GH}((M,p),(N,q))<\epsilon$, we fix a metric $d$ on the disjoint union space $M\amalg N$ such that $d((M,p),(N,q))<\epsilon$. For each $x\in M$ one can define a map $h:M\to N$ such that $h(x)$ is a point in N with $d(h(x),x)<\epsilon$. Thus the triangle inequality gives that $h$ is $4\epsilon$-continuous (cf. {\bf [P1]}). Hence by [Main obstruction result 3 in {\bf [P1]}] there is a continuous map $f:M\to N$ with $d(h(x),f(x))\leq (n+2)\epsilon$ for all $x\in M$. Take $\tau(\epsilon)=4(n+3)\epsilon$. It is not difficult to check that $f:M\to N$ is an  $\tau(\epsilon)$-Hausdorff approximation with $\lim_{\epsilon\to 0}\tau(\epsilon)=0$. Moreover, from [Corollary 4.6 in {\bf [P1]}] and the second-named author's argument in {\bf [W]} we have the induced map $f_*:\pi_1(M,p)\to\pi_1(N,q)$ is surjective and hence (3.1.3) is established.$\Box$\\
 
The proof of (3.1.4) is basically along the line of Fukaya and Yamaguchi's proof in [Theorem 7.1 in {\bf [FY1]}]. Here we will point out how the original process can work under our settings. For the sake of brief, we say a group has {\bf property (*)} if there is a subgroup satisfying (i) and (ii) in (3.1.4).\\
\\
{\bf Proof of (3.1.4).} This proof is done by induction on $\dim N$ and by contradiction. When $\dim N=n$, by [Theorem A1.12 in {\bf [CC2]}], we can choose $\epsilon$ small enough such that there is a diffeomorphism from $(M^n,p)$ to $(N^n,q)$, which we take as $f$. Since $\pi_1(V)$ is trivial, $\Gamma$ is also trivial and the theorem holds in this case. Now we suppose (3.1.4) holds for $k<\dim N<n$ with fixed $k$ but not hold for $\dim N=k$. Then for sequences $\epsilon_i\to 0$ and $w_i\to\infty$ as $i\to\infty$, there exists sequences $(M_i^{n},p_i)$ and $(N_i^{k},q_i)$ satisfying (3.1.1) and (3.1.2) but no map $(M_i,p_i)\to (N_i,q_i)$ satisfies (3.1.4) for $\epsilon=\epsilon_i$ and $w=w_i$ simultaneously. Note that when $i$ large enough we always have a continuous $\tau(\epsilon_i)$-Hausdorff approximation map $f_i: M_i\to N_i$ with $f_i(p_i)=q_i$.

 Let $V_i=B_{q_i}(\frac{\mu_0}{2})$, $U_i=f_i^{-1}(V_i)$ and $\Gamma_i=\pi_1(U_i)$ be defined as above.  In order to use the 
induction hypothesis, we need to blow-up the metrics as the technique shown in {\bf [FY1]}. However, since $f_i$ 
is only a continuous map in our case, we consider a scaling of metrics as follows. Since $\dim(M_i^n)>\dim(N_i^k)$,
 there exists $q_i'\in B_{q_i}(\frac{\tau(\epsilon_i)}{10})$ such that $diam(f_i^{-1}(q_i'))>0$. Let $p_i'\in f_i^{-1}(q_i')$. Then it 
can be shown that 
\[d_H((M_i^n,p_i'),(N_i^k,q_i'))<2\tau(\epsilon_i).\]

Indeed, since $d_H((M_i^n,p_i),(N_i^k,q_i))<\epsilon$, there exist $\epsilon_i$-pointed Hausdorff approximations $h_i$ and $\tilde{h}_i$ such that
the map $f_i$ is induced by $h_i$ as in the proof of (3.1.3) and $h_i(p_i)=q_i$, $\tilde{h}_i(q_i)=p_i$, 
$h_i(B_{p_i}(\frac{1}{\epsilon_i}))\subseteq B_{q_i}(\frac{1}{\epsilon_i}+\epsilon_i)$ and
 $\tilde{h}_i(B_{q_i}(\frac{1}{\epsilon_i}))\subseteq B_{p_i}(\frac{1}{\epsilon_i}+\epsilon_i)$. 
Moreover, $d(p_i,p_i')<\frac{6}{5}\tau(\epsilon_i)$
since $f_i$ is a $\tau(\epsilon_i)$-Hausdorff approximation. Note that $\epsilon_i<\frac{\tau(\epsilon)}{10}$ and then 
for $x\in B_{p_i'}(\frac{1}{2\tau(\epsilon_i)})$ we have
\begin{eqnarray*}
d(h_i(x),q_i')&\leq&d(h_i(x),h_i(p_i))+d(q_i,q_i')\\
              &<&d(x,p_i)+2\epsilon_i+d(q_i,q_i')\\
              &\leq&d(x,p_i)+d(p_i,p_i')+d(q_i,q_i')+2\epsilon_i\\
              &<&\frac{1}{2\tau(\epsilon_i)}+2\tau(\epsilon_i)
\end{eqnarray*}
Thus $h_i(B_{p_i'}(\frac{1}{2\tau(\epsilon_i)}))\subseteq B_{q_i'}(\frac{1}{2\tau(\epsilon_i)}+2\tau(\epsilon_i))$. Similarly, 
we can show that $\tilde{h}_i(B_{q_i'}(\frac{1}{2\tau(\epsilon_i)}))\subseteq B_{p_i'}(\frac{1}{2\tau(\epsilon_i)}+2\tau(\epsilon_i))$
and then $d_H((M_i^n,p_i'),(N_i^k,q_i'))<2\tau(\epsilon_i).$ Therefore, for the sake of brief, we can assume the map
$f_i:M_i^n\to N_i^k$ with $f_i(p_i)=q_i$ satisfying the condition as in (3.1.3)and $diam(f_i^{-1}(q_i))>0$ for $i$ large enough.

For each $i$ denote $\delta_i=\sup\{d(x,y)|x,y\in f_i^{-1}(q_i)\}$ Then $\delta_i\to 0$ as $i\to\infty$. 
Blow-up the original metrics $g_{M_i}$ and $g_{N_i}$ of $M_i$ and $N_i$ respectively by 
\[g_i=\frac{1}{\delta_i}g_{M_i} \;\;\;and\;\;\; h_i=\frac{1}{\delta_i}g_{N_i}.\] 
Then, by taking a subsequence if necessary, $((V_i,h_i),q_i)$ converges to $(\mbox{\boldmath$R$}^k,can,0)$ since $Ric_{N_i}\geq -(n-1)$ and $inj_{q_i}(N_i)>\mu_0$. Assume $((U_i,g_i),p_i)$ converges to a pointed metric space $(X,x_0)$. Moreover, we have that the sequence $\{f_i\}$ is an {\em almost equicontinuous family} and convergent.\\
\\
{\bf Sublemma 3.3.} {\em $f_i$ converges to a continuous map $f:X\to\mbox{\boldmath$R$}^k$ with}\\
             (3.3.1) {\em For every $x,y\in\mbox{\boldmath$R$}^k$, $d(f^{-1}(x),f^{-1}(y))=d(x,y)$.}\\
             (3.3.2) {\em For every $x'\in f^{-1}(x)$ there exists a point $y'\in f^{-1}(y)$ such that                                             $\;\;\;d(x',y')=d(x,y)$.}\\
             (3.3.3) {\em $f(x_0)=0$.}\\
\\
{\bf Proof of Sublemma 3.3.} We may assume that $\tau(\epsilon_i)$ decreases monotonically to $0$. Fix $i_0$ and choose
            $\eta_{i_0}=\frac{3}{2}\tau(\epsilon_{i_0})$. Then given $\eta >\eta_{i_0}$ we have
            $d(f_i(x_i),f_i(y_i))<\eta$ for all $x_i,y_i\in U_i$
            provided $d(x_i,y_i)<\frac{1}{2}\tau(\epsilon_i)$ and $i\geq i_0$. Take a dense subset
            $A_i=\{a_1^{i},a_2^{i},\ldots\}\subset U_i$ with $a_j^{i}\to a_j\in X$. Then the set $A=\{a_j\}_{j=1}^{\infty}\subset X$ is
            also dense in $X$. By use of the diagonal process, the map $f:A\to\mbox{\boldmath$R$}^k$ defined by
            $f(a_j)\equiv\lim_{i\to\infty}f_i(a_j^{i})$ is well-defined.
            Moreover, by taking $i_0\to\infty$, we have $\{f(a_j)\}$ is a Cauchy sequence and hence $f$ is
            uniformly continuous from $X$ to $\mbox{\boldmath$R$}^k$. Clearly $f$ satisfies equations (3.3.1)-(3.3.3).$\Box$\\ 

            From Theorem 2.4 and (3.3.1)-(3.3.3), we apply [Lemma 7.4 and 7.5 in {\bf [FY1]}] to conclude that $X$ is isometric to a product $\mbox{\boldmath$R$}^k\times Z$, where $Z$ is compact and not a single point, and the map $f:\mbox{\boldmath$R$}^k\times Z\to \mbox{\boldmath$R$}^k$ is in fact the projection.

            Let $B^{k}(r_0)$ denote the metric $r_0$-ball in $\mbox{\boldmath$R$}^k$ around the origin. Denote  $U_i(r_0)\equiv f_i^{-1}(B_{q_i}(r_0))$, where $B_{q_i}(r_0)\subseteq (N_i,h_i)$. Then                        
 $$\lim_{i\to\infty}U_i(r_0)=B^{k}(r_0)\times Z.\leqno{(3.4)}$$ 
Let $d_0$ be the distance of $B^{k}(r_0)\times Z$. By Corollary 2.6, one can find sequences $y_j\in B^{k}(r_0)\times Z$ and $r_j\to\infty$ as $j\to\infty$ such that
$$\lim_{j\to\infty}((B^{k}(r_0)\times Z,r_jd_0),y_j)\;=\;((\mbox{\boldmath$R$}^{m},can),0))\leqno{(3.5)}$$
, where $m>k$ since $Z$ is not one point. Now combining (3.4) and (3.5) together, it is not difficult to conclude that for given $\epsilon>0$ there are $i_0$, $j_0$ and $\hat{p}_i\in U_i(r_0)$ such that for $i\ge i_0$ we have
 $$d_H((U_i(r_0),r_{j_0}g_i),\hat{p}_i),((\mbox{\boldmath$R$}^{m},can),0))\;<\;\epsilon.\leqno{(3.6)}$$ 
Therefore by induction hypothesis, we have, for $i$ large enough, a map $\Phi_i:(U_i(r_0),r_{j_0}g_i),\hat{p}_i)\to((\mbox{\boldmath$R$}^{m},can),0)$ such that $\Phi_i$ is a $\tau(\epsilon)$-Hausdorff approximation and the fundamental group $\pi_1(\Phi_i^{-1}(B^m(\frac{\mu_0}{2})))$ satisfying property (*) for $w=w_{n,m}$.

Let $\Gamma_i(\hat{p}_i,\mu_0)=Im[i_*:\pi_1(\Phi_i^{-1}(B^m(\frac{\mu_0}{2})))\to \pi_1(U_i)]$ be the image of the induced map $i_*$ of the inclusion map $i:\Phi_i^{-1}(B^m(\frac{\mu_0}{2}))\to U_i$. Then  $\Gamma_i(\hat{p}_i,\mu_0)$ has naturally the property (*).

Let $(\tilde{U_i},\tilde{g_i},\tilde{p_i})$ be the universal Riemannian covering space of $(U_i,g_i,p_i)$ with covering map $\Pi_i:(\tilde{U_i},\tilde{g_i},\tilde{p_i})\to (U_i,g_i,p_i)$. Then $\Gamma_i=\pi_1(U_i,p_i)$ is the deck transformation group. By taking a subsequence if necessary, we may assume that there exists a triple $(W,G,\tilde{p}_{\infty})\in {\cal M}_{eq}$ such that 
$$\lim_{i\to\infty}(\tilde{U_i},\Gamma_i,\tilde{p_i})\;=\;(W,G,\tilde{p}_{\infty}),\leqno{(3.7)}$$
$$\Pi_i\; converges\; to\; a\; map\; \Pi_{\infty}:W\to \mbox{\boldmath$R$}^{k}\times Z.\leqno{(3.8)}$$
Since $\Pi_{\infty}$ also satisfies (3.3.1) and (3.3.2), by Theorem 2.4 and Corollary 2.5 we have $W$ is isometric to $\mbox{\boldmath$R$}^{k}\times W'$ and $W'$ is isometric to $(\mbox{\boldmath$R$}^{\ell}\times Y)$, where $Y$ is a compact metric space.

Let $\tilde{U_i}(r_0)=\Pi_i^{-1}(U_i(r_0))$. Then
$$\lim_{i\to\infty}(\tilde{U_i}(r_0),\Gamma_i,\tilde{p_i})\;=\;(B^k(r_0)\times\mbox{\boldmath$R$}^{\ell}\times Y,G,\tilde{p}_{\infty})\leqno{(3.9)}$$

Apply Theorem 2.8, we have for each $\varepsilon>0$ a normal subgroup $G_{\varepsilon}$ of $G$ such that (2.8.1)-(2.8.3) hold. Therefore, Theorem 2.3 gives a sequence of normal subgroups $\Gamma_{i,\varepsilon}$ of $\Gamma_i$ with
$$\lim_{i\to\infty}(\tilde{U_i}(r_0),\Gamma_{i,\varepsilon},\tilde{p_i})\;=\;(B^k(r_0)\times\mbox{\boldmath$R$}^{\ell}\times Y,G_{\varepsilon},\tilde{p}_{\infty}),\;and\leqno{(3.10)}$$
$$\Gamma_i/\Gamma_{i,\varepsilon}\cong G/G_{\varepsilon}\leqno{(3.11)}$$ 
for each sufficiently large $i$.

Now we investigate the relationship between $\Gamma_{i,\varepsilon}$ and $\Gamma_i(x,\varepsilon)$ for some $\varepsilon>0$. \\
\\ 
{\bf Sublemma 3.12.} {\em For every $x\in U_i(r_0)$, $\Gamma_{i,\varepsilon}\subset\Gamma_i(x,3\varepsilon)$ for $i$ large enough.}
\\
{\bf Proof of Sublemma 3.12.} Let $\tau_i$ be the equivarient pointed Hausdorff distance between $(\tilde{U_i}(r_0),\Gamma_{i,\varepsilon},\tilde{p_i})$ and $(B^k(r_0)\times\mbox{\boldmath$R$}^{\ell}\times Y,G_{\varepsilon},q)$. To prove this lemma, we shall show that, for $i$ sufficiently large, each element $\gamma\in\Gamma_{i,\varepsilon}$ can be generated by geodesic loops at $x$ of length less than $C\tau_i+\epsilon$, where $C=C(G,\varepsilon)>0$. 

As in Definition 2.1, let $\varphi_i:(\tilde{U_i}(r_0),\Gamma_{i,\varepsilon},\tilde{p_i})\to (B^k(r_0)\times\mbox{\boldmath$R$}^{\ell}\times Y,G_{\varepsilon},q)$ be a $\tau_i$-Hausdorff approximation and let $\lambda_i:\Gamma_{i,\varepsilon}(\frac{1}{\tau_i})\to G_{\varepsilon}(\frac{1}{\tau_i})$, $\lambda_i':G_{\varepsilon}(\frac{1}{\tau_i})\to \Gamma_{i,\varepsilon}(\frac{1}{\tau_i})$ be the corresponding maps. Note that we can take a point $\tilde{x}\in\tilde{M_i}$ over $x$ such that $d(\tilde{x},\tilde{p_i})$ is uniformly bounded, say $2r_0$. By (2.3.8) we may assume the length of $\gamma$ is uniformly bounded by a constant. Hence we can assume $\bar{\gamma}=\lambda_i(\gamma)$ for each sufficiently large $i$. By (2.8.3) there are $\bar{\gamma_1},\ldots,\bar{\gamma_s}\in G_{\varepsilon}$ such that
$$\bar{\gamma}=\bar{\gamma_s}\bar{\gamma}_{s-1}\ldots\bar{\gamma_1},\leqno{(3.12.1)}$$
$$d(\bar{\gamma_j}\bar{\gamma}_{j-1}\ldots\bar{\gamma_1}(\varphi(\tilde{x})),\bar{\gamma}_{j-1}\ldots\bar{\gamma_1}(\varphi(\tilde{x}))<\varepsilon, \;for\; 1\leq j\leq s.\leqno{(3.12.2)}$$
For each $j$ we put $\gamma_j=\lambda_i'(\bar{\gamma_j})$. Then $\gamma$ has the expression
$\gamma\;=\;\gamma_s\ldots\gamma_1(\gamma_s\ldots\gamma_1)^{-1}\gamma.$ and therefore,
$$d(\gamma_j\gamma_{j-1}\ldots\gamma_1(\tilde{x}),\gamma_{j-1}\ldots\gamma_1(\tilde{x}))<2(j+1)\tau_i+\varepsilon\leqno{(3.12.3)}$$
$$d(\gamma(\tilde{x}),\gamma_s\ldots\gamma_1(\tilde{x})<2(s+1)\tau_i.\leqno{(3.12.4)}$$
Note that $s$ depends on $G$, $\varepsilon$ and $\mu_0$. Then the sublemma holds.$\Box$\\

Up to now we know that $\Gamma_{i,\mu_0}$ has a subgroup $H_{i,\mu_0}$ with the property (*). Then we can apply the argument in [{\bf [FY1]}, p.288 and p.289] to conclude that $\pi_1(U_i)$ has a subgroup $H$ with property (*) for $i$ large enough. It contradicts to our assumption. Hence we have, for $\dim N=k$, there exist $\epsilon_{n,k}(\mu_0)$, $w_{n,k}$ and a map $f_i:{M_i}^n\to {N_i}^k$ such that (3.1.3) and (3.1.4) hold under the assumption of (3.1.1) and (3.1.2). Thus, by induction we finish the proof of the Technical Lemma 3.1.$\Box$\\
   \section{A weak Margulis's Lemma}      

In this section we will give a proof of Theorem B. First, by Technical lemma 3.1, we have the following lemma.\\
\\
{\bf Lemma 4.1.} {\em For given integer $n$ and $k$ with $n\ge k$ there exist $\delta_{n,k}$, $I_{n,k}$ and $w_{n,k}$ depending only on $n$ and $k$ satisfying the following: Let $(M_i^n,p_i)$ be a sequence of complete pointed Riemannian $n$-manifold with $Ric_{M_i^n}\ge -(n-1)$. Suppose $(M_i^n,p_i)$ converges to a metric space $(X,p_{\infty})$ of dimension $k$ in the pointed Hausdorff convergence. Then there exist $p_i'\in B_{p_i}(1/2)$ such that, for $\delta<\delta_{n,k}$ and $i\ge I_{n,k}$, the image of the inclusion homomorphism
\[\Gamma'\;=\;Im[\pi_1(B_{p_i'}(\delta))\,\to\,\pi_1(B_{p_i}(1))]\]
admits a subgroup $\Lambda'\subset\Gamma'$ with}\\
(4.1.1) {\em $\;\;[\Gamma':\Lambda']<w_{n,k}$};\\
(4.1.2) {\em $\;\;\Lambda'$ is solvable with length of polycyclicity$\leq n-k$}.
\\
{\bf Proof of Lemma 4.1.} We prove it by contradiction. Assume that there exit $\delta_i\to 0$, $w_i\to\infty$ and $M_i$ with $Ric_{M_i}\ge -(n-1)$ satisfying the following: For each $p_i'\in B_{p_i}(1/2)$ and each $I>0$, there exists an $i\ge I$ such that 
\[\Gamma_{i}'\;=\;Im[\pi_1(B_{p_i'}(\delta_{i}))\,\to\,\pi_1(B_{p_i}(1))]\]
never admit a subgroup with properties (4.1.1) and (4.1.2) in Lemma 4.1 for $w=w_i$.

By Corollary 2.6 there are sequences $q_i'\in B_{p_i}(1/2)$  and $r_i\to\infty$ such that $\lim_{i\to\infty}((M_i,r_ig_i),q_i')=(\mbox{\boldmath$R$}^{m},can,0)$, where $m\ge k=\dim(X)$.
Note that the sequence $r_i$ can be chosen so that $r_i\delta_{i}\to 0$ as $i\to\infty$. By Technical lemma 3.1, there exist an $i_0$ large enough and a Hausdorff approximation $f_{i_0}:(M_{i_0}^n,q_{i_0}')\to (\mbox{\boldmath$R$}^{k},0)$ with $f_{i_0}(q_{i_0}')=0$ so that the fundamental group $\Gamma_{i_0}''=\pi_1(f_{i_0}^{-1}(B^k(10)))$ admits a subgroup $\Lambda_{i_0}''\subset\Gamma_{i_0}''$ satisfying properties (4.1.1) and (4.1.2) of Lemma 4.1 for $w_{n,k}$ independent of $i$. By $r_i\delta_{i}\to 0$, one has
\[\Gamma_{i_0}'\;\subset\;Im[\Gamma_{i_0}''\to\pi_1(M_{i_0})].\]
Then $\Gamma_{i_0}'$ admits a subgroup $\Lambda_{i_0}'$ satisfying properties (4.1.1) and (4.1.2) of Lemma 4.1. This is a contradiction.$\Box$\\
\\
{\bf Proof of Theorem B.} Now we give a proof of Theorem B. For given a divergent sequence $r_i\to\infty$ there exists a metric space $(X,p_{\infty})$ of dimension $k\leq n$ such that $((M^n,r_ig_M),p)$ converges to $(X,p_{\infty})$. By Lemma 4.1 choose $\delta_n=\min_{0\leq k\leq n}{\delta_{n,k}}$ and $w_n=\max_{0\leq k\leq n}{w_{n,k}}$ and then Theorem B holds.$\Box$\\
\\
{\bf Remark 4.2.} Consider a compact Riemannian $n$-manifold $M^n$ with $Ric_{M^n}\ge 0$. Scaling the metric of $M$ so that $diam(M)\leq\frac{\delta_n}{2}$ and we still have $Ric_M\ge 0$. Then we can conclude that the fundamental group $\pi_1(M)$ of $M$ admits a subgroup $H$ such that \\
(4.2.1) $\;\;\;H$ is solvable with length of polycyclicity $\leq n$;\\
(4.2.2) $\;\;\;[H:\pi_1(M)]\;<\;w_n$.\\
 In the next section, we will use this result to give a proof of Theorem A.\\
    \section{Proof of Theorem A}       
Now we are in a position to prove Theorem A. We divide the proof into the proof of part (a) and part (b) and follow  similar methods in the proofs of [Corollary 7.20 in {\bf [FY1]}, p.289] and [Corollary 0.11 in {\bf [FY1]}, p.290 and p.291] respectively. However, we extend their results to manifolds with positive Ricci curvature.\\
\\
{\bf Proof of (a).} Let $M^n$ be a compact Riemannian $n$-manifold with $Ric\ge 0$. By Remark 4.2, the fundamental group $\pi_1(M)$ of $M$ admits a subgroup $H$ satisfying (4.2.1) and (4.2.2) and then $b_1(M,\mbox{\boldmath$Z$}_p)\leq n$ for all prime $p\geq w_n\equiv p_n$. Here we give a proof that the maximal case can never occur. 

Our proof is by contradiction. We use the same notations as in the proof of Technical Lemma 3.1 and consider the case when $k=0$. Then we have $U_i=M_i$, $\tilde{U_i}=\tilde{M_i}$ and         
 \[(\mbox{\boldmath$R$}^{\ell}\times Y)/G\;\,=\;\,Z.\]
by (3.4), (3.7), (3.8) and (3.9). Therefore, $\mbox{\boldmath$R$}^{\ell}/\phi(G)$ is a compact, where the projection $\phi:G\to Isom(\mbox{\boldmath$R$}^{\ell})$ is defined as in Theorem 2.8. By [Corollary 4.2 in {\bf [FY1]}, p.273], the generalized Bieberbach's theorem, there is a finite-index normal subgroup $G^{'}$ of $\phi(G)$ such that $\mbox{\boldmath$R$}^{\ell}/G^{'}$ is a flat $s$-torus $T^{s}$. Moreover, by Theorem 2.3, we put a finite-index normal subgroup $\hat{\Gamma}_i$ of $\Gamma_i=\pi_1(M)$ converging to $\phi^{-1}(G^{'})$. Thus the compact manifold $\hat{M}_i\equiv\tilde{M}_i/\hat{\Gamma}_i$ converges to $T^{s}$. Note that $Ric_{\hat{M}_i}>0$ and $\hat{M}_i$ is also a covering space of $M$. Then $\pi_1(\hat{M_i})$ is a subgroup of $\pi_1(M_i)$. Use (3.1.3) we have a surjection between $\pi_1(\hat{M_i})$ and $\pi_1(T^{s})$ for $i$ large enough. Thus $\pi_1(\hat{M_i})$ has infinite order and it contradicts to $\pi_1(M)$ is finite. Hence $b_1(M,\mbox{\boldmath$Z$}_p)\leq n-1$ for all prime $p\geq w_n\equiv p_n$.$\Box$\\
\\
{\bf Proof of (b).} We also prove it by contradiction. Suppose there exists a sequence of Riemannian manifolds $M_i$ 
        with $Ric_{M_i}>0$, and universal Riemannian covering $\tilde{M_i}$ with, 
        \[\lim_{i\to\infty}\frac{diam(\tilde{M_i})}{diam(M_i)}\;=\;\infty.\]    
By scaling the metric, we may assume $diam(\tilde{M_i})=1$ for all $i$. Then there is a sequence $\epsilon_i$ with $\lim_{i\to\infty}\epsilon_i=0$ such that $diam(M_i)<\epsilon_i$. Note that for each $i$ $\pi_1(M_i)$ admits a subgroup $H_i$ satisfying (4.2.1) and (4.2.2).

  Define $H_i^{0}=H_i$, $H_i^{j}=[H_i^{j-1},H_i^{j-1}]$, the commutator of $H_i^{j-1}$, and then $H_i^{n}=0$. Choose $\tilde{p_i}\in\Pi_{i}^{-1}(p_i)$. Since $diam(\tilde{M_i}/ H_i^{j})\leq 1$ for all $i$ and $j$, Corollary 2.5 implies there exists a triple $(X,G_j,x_0)$ such that             
\[\lim_{i\to\infty}(\tilde{M_i},H_i^{j},\tilde{p_i})\,=\,(X,G_j,x_0)\]            
 for each $0\leq j\leq n$. Note that $G_0$ acts on $X$ transitively. Take a number $j_0$ such that $X/G_{j_0-1}$ is a point and $X/G_{j_0}$ is not a point.

 Now we put 
\[\bar{M_i}\,=\,\tilde{M_i}/H_i^{j_0},\;\;\Lambda_i\,=\,H_i^{j_0}/H_i^{j_0-1}.\]            
Let $\bar{p_i}\in\bar{M_i}$ be the point corresponding to $\tilde{p_i}$. By Corollary 2.6 and repeating the blow-up arguments (at most finite times), we may assume $(\bar{M_i},\bar{p_i})$ converges to $(\mbox{\boldmath$R$}^{\ell},0)$ for some $\ell>0$. Let $\Lambda$ be the group such that $\lim_{i\to\infty}(\bar{M_i},\Lambda_i,\bar{p_i})=(\mbox{\boldmath$R$}^{\ell},\Lambda,0).$ Then the abelian group $\Lambda$ acts on $\mbox{\boldmath$R$}^{\ell}$ transitively and, in fact, $\Lambda$ is just the vector group $\mbox{\boldmath$R$}^{\ell}$.
                   
Now we also use the {\em pseudogroup technique} as in {\bf [Y]}. Let 
\[\Lambda_{i}'=\{\gamma\in\Lambda_i\mid d(\gamma(\bar{p_i}),\bar{p_i})<10\ell\}.\]
 Consider $B^{\ell}(10\ell)\subset\mbox{\boldmath$R$}^{\ell}$ and $\Lambda_{i}'$ as the pseudogroups of isometric embeddings of $B^{\ell}(10\ell)$ to $B^{\ell}(20\ell)$ and $B_{\bar{p_i}}(10\ell)$ to $B_{\bar{p_i}}(20\ell)$ in $\bar{M_i}$ respectively. Consider the lattice
\[E_{\infty}\;=\;B^{\ell}(10\ell)\cap\mbox{\boldmath$Z$}^{\ell}\]
and take $\gamma_{1_i},\ldots,\gamma_{\ell_i}\in\Lambda_{i}'$ such that $\gamma_{j_i}$ converges to $e_j\in E_{\infty}$, where $e_1,\ldots,e_{\ell}$ are the canonical basis of $\mbox{\boldmath$Z$}^{\ell}$. We denote by $E_i$ the pseudogroup of $\Lambda_{i}'$ generated by $\gamma_{1_i},\ldots,\gamma_{\ell_i}$. Since $E_i$ is abelian, $(B_{\bar{p_i}}(10\ell),E_i)$ converges to $(B^{\ell}(10\ell),E_{\infty})$. It follows that $\hat{M_i}\equiv  B_{\bar{p_i}}(10\ell)/E_i$ converges to the flat torus $T^{\ell}=B^{\ell}(10\ell)/E_{\infty}$ with respect to the Hausdorff distance. Note that, as in the proof of part (a), $Ric_{\hat{M}_i}>0$ and $\hat{M_i}$ is also a covering of $M_i$. Again (3.1.3), there is a surjection from $\pi_1(\hat{M_i})$ to $\pi_1(T^{\ell})$ for $i$ large enough and then the group $\pi_1(\hat{M_i})$ has infinite order. Thus $\pi_1(M_i)$ also has  infinite order and this is a contradiction. Therefore we complete the proof.$\Box$\\

\end{document}